\documentclass[12pt, reqno]{amsart} 
\usepackage{amssymb,amscd,amsfonts,amsbsy}
\usepackage{latexsym}
\usepackage{exscale}
\usepackage{amsmath,amsthm,amsfonts}
\usepackage{mathrsfs}
\usepackage{xcolor} 
\usepackage[colorlinks=true,linkcolor=blue,citecolor=red,urlcolor=red, backref=page]{hyperref} 
\usepackage{esint} 

\usepackage[utf8]{inputenc}

\calclayout

\usepackage{verbatim} 
\usepackage{esint}

\numberwithin{equation}{section}

\DeclareMathOperator{\hdim}{dim_H}

\usepackage{xcolor}

\numberwithin{equation}{section}

\theoremstyle{plain}

\newtheorem{theorem}{Theorem}[section]
\newtheorem{lemma}[theorem]{Lemma}
\newtheorem{proposition}[theorem]{Proposition}
\newtheorem{corollary}[theorem]{Corollary}

\theoremstyle{definition}

\newtheorem{remark}[theorem]{Remark}

\theoremstyle{remark}
\numberwithin{equation}{section}

\newcommand{\e}{\varepsilon}

\newcommand{\supp}{\textrm{supp}\ }

\newcommand{\re}{\textrm{Re}}
\newcommand{\im}{\textrm{Im}}

\newcommand{\R}{\mathbb{R}}
\newcommand{\Z}{\mathbb{Z}}
\newcommand{\C}{\mathbb{C}}

\newcommand{\eps}{\varepsilon}
\newcommand{\norm}[1]{\left\lVert #1 \right\rVert}

\def\XXint#1#2#3{{\setbox0=\hbox{$#1{#2#3}{\int}$}
    \vcenter{\hbox{$#2#3$}}\kern-.5\wd0}}

\newcommand{\vertiii}[1]{{\left\vert\kern-0.25ex\left\vert\kern-0.25ex\left\vert #1 
    \right\vert\kern-0.25ex\right\vert\kern-0.25ex\right\vert}}

\numberwithin{equation}{section}

\author{Carolina A. Mosquera}
\author{Andrea Olivo}
\address{
Departamento de Matem\'atica, Facultad de Ciencias Exactas y Naturales, Universidad de Buenos Aires \\
Ciudad Universitaria, Pabell\'{o}n I (C1428EGA) \\
 Ciudad de Buenos Aires, Argentina and IMAS-CONICET, Argentina}
\email{mosquera@dm.uba.ar}
\address{
       ICTP - The Abdus Salam International Centre for Theoretical Physics \\
       Strada Costiera, 11, I - 34151 \\
       Trieste, Italy}
 \email{aolivo@ictp.it}
\thanks{CM is partially supported by grants PICT 2018--03399 and PICT 2018--04027.}
\subjclass[2010]{28A78, 28A80}
\keywords{Fourier decay, homogeneous self-similar measures, correlation dimension}

\begin{document}

\title{Fourier decay of self-similar measures on the complex plane}
\maketitle

\begin{abstract}
We prove that the Fourier transform of self-similar measures on the complex plane has fast decay outside of a very sparse set of frequencies, with quantitative estimates, extending the results obtained in the real line, first by R. Kaufman, and later, with quantitative bounds, by the first author and P. Shmerkin.
Also we derive several applications concerning correlation dimension and Frostman exponent of these measures. Furthermore, we present a generalization for a particular case on $\R^n,$ with $n\ge3.$

\end{abstract}

\section{introduction and statement of results}
Given $\mu$ a finite Borel measure on the complex plane, its Fourier transform is defined as
\[
\widehat{\mu}(\xi)=\int_{\C} e^{2\pi i \re(z\cdot\overline{\xi})}\,  d\mu(z).
\]
The behaviour of $\widehat{\mu}(\xi)$ when $|\xi| \to \infty$ is a fundamental characteristic of the measure $\mu$. Measures for which $|\widehat{\mu}(\xi)| \to 0$ when $|\xi|\to \infty$ are called \textit{Rachjman measures}. By the Riemann-Lebesgue lemma every absolutely continuous measure is Rachjman, but also many singular measures. Among the (possibly) singular measures, an important group are the homogeneous self-similar measures (see Section \ref{Decay}), which are the measures that we consider in this paper. In particular, we focus on self-similar measures supported in the complex plane.

The decay properties of $\widehat{\mu}(\xi)$ as $|\xi|\to +\infty$ give a crucial ``arithmetic" information about $\widehat{\mu}(\xi).$ 
We say that $\widehat{\mu}(\xi)$ has polynomial decay if there exist $\sigma, C_{\sigma}>0$ such that
$
|\widehat{\mu}(\xi)|\le C_{\sigma} |\xi|^{-\sigma/2}.
$
Even more, the \textit{Fourier dimension} of $\mu$ is defined as
\begin{equation}
\label{eq:power-decay}
\dim_{F}(\mu) =   2\sup\{  \sigma\ge 0: |\widehat{\mu}(\xi)| \le C_\sigma |\xi|^{-\sigma} \text{ for some $C_\sigma>0$ and all $\xi\neq 0$ }\},
\end{equation}
and then one can say that  $\mu$ has polynomial decay if and only if $\dim_{F}(\mu)>0$.
For many purposes simple convergence of $\widehat{\mu}$ to zero is not enough, and some quantitative decay is needed. For example, if $\dim_{F}>0$, then
$\mu$-almost every number is normal to any base, see \cite{DavenportEtAl63, QR03}, and if $\dim_{F}(\mu) > 0$
and $\mu$ satifies a Frotsman type condition, then $\mu$ satisifes a restriction theorem analog to the Stein-Tomas theorem for the sphere,  see \cite{Mitsis02,Mockenhaupt00}.

Regardless of its significance, the Fourier dimension of a measure is particularly difficult to calculate or even give some estimative. Furthermore, in some cases there are not decay at all. For example, consider the complex Bernoulli convolutions $\mu_{\lambda}$, being the distribution of the random series 
$\sum_{n=1}^{+\infty}  \pm \lambda^n,$ where the signs are chosen independently with probabilities $\{1/2, 1/2\}$ and $\lambda \in \mathbb{D}$, the open unit disk.
Also, $\mu_{\lambda}$ can be defined as the self-similar measure associated to the iterated function system (IFS) $\{\lambda z-1, \lambda z +1\}.$
 Solomyak and Xu \cite{SX03} proved that if $\theta$ is a complex Pisot number and $1<|\theta|<\sqrt{2},$ then 
$|\widehat{\mu}_{\lambda}(\xi)|\nrightarrow 0$ as $|\xi| \to\infty$, which implies that $\dim_F \mu_\lambda =0$ (and in particular, $\mu_{\lambda}$ is singular for $\lambda=1/{\theta}$).
Recall that a non-real algebraic integer $\theta$, with $|\theta|>1$, is called a complex Pisot
number if all its Galois conjugates, except $\theta$, are less than one in modulus.
Later, Shmerkin and Solomyak \cite{SS16} proved that the Fourier transform of complex Bernoulli convolutions have power decay for all parameter $\lambda$, outside of an exceptional set of
parameters of zero Hausdorff dimension.

Nevertheless, if a measure $\mu$ has zero Fourier dimension, it may happen that $\widehat{\mu}$ has fast decay outside of a very sparse set of frequencies. In fact, Kaufman \cite{Kaufman84}, for homogeneous self-similar measures on the real line and Tsujii \cite{Tsujii15}, for the non-homogeneous case,
proved that for any $\varepsilon> 0$ there exists
$\delta >0$ such that the set
\[
\{\xi \in [-T, T]\colon\widehat{\mu}(\xi) \geq T^{-\delta} \}
\]
can be covered by $T^{\varepsilon}$ intervals of side length 1. Kaufman use a version of the well known Erd\"{o}s-Kahane argument whereas the proof of Tsujii is based on large deviation estimates.
 Recently, the first author and Shmerkin \cite{MS18} made the dependance of $\delta$ on $\varepsilon$ quantitative in the homogeneous case.

Even more, if there is not decay at all, in \cite{Kaufman84}, for the case of Bernoulli convolutions in the real line with parameter $\lambda \in (0, 1/2)$ and then, generalized to any non-atomic homogeneous self-similar measure in the real line in \cite{MS18}, the authors proved that, if $F''>0$ is a $C^2$ diffeomorphism of $\mathbb{R}$ with $F''>0$, then $\dim_F{F\mu}>0$, where here and below $F\mu (A)= \mu(F^{-1}(A))$ for all borel sets $A \subset \mathbb{R}$ denotes the push-forward measure.

 The goal of this paper is to extend the  results mentioned above to the case of homogeneous self-similar measures on the complex plane. On one hand, we prove that the Fourier transform of $\mu$ has fast decay outside of a very sparse set of frequencies for homogeneous self-similar measures on the complex plane, extending the results obtained in \cite{MS18} for the real case, making the dependence of $\delta$ on $\varepsilon$ explicit (see Proposition \ref{prop:EK0} and Proposition \ref{prop:EK0_real}). Although the approach is similar to the one used in \cite{MS18} and \cite{Kaufman84} and is based on the Erd\"{o}s-Kahane argument, but there are some new features, mainly
because we have to deal with two different cases, when the complex contraction ratio of the iterated function system belongs or not to the real line.

Then, as an application, we obtain a generalization of the Kaufman's result. We prove that if $F$ is an analytic function with $F''\neq 0$ in a neighborhood of the $\supp \mu$, then $F\mu$ has polynomial decay with quantitative estimates (see Theorem \ref{Kaufmann on C}).

 Given $\nu$ a Borel measure, since convolution (of a measure, of a function, etc) with $\nu$ is a smoothing operation, a natural problem is quantifying the additional degree of smoothness ensured by convolving with $\nu.$ In \cite{RS20} the authors prove that uniformly perfect measures on the real line (which includes Ahlfors-regular measures as a proper subset) have the property that convolving with them results in a strict increase of the $L^q$ dimension. 
 In \cite{MS18} the authors consider the special case when $\nu$ is a homogeneous self-similar measure on the real line and $q=2,$ with quantitative estimates. Also non-quantitative results can be also deduced before from \cite{Shmerkin16}.

In the present work we prove that convolving with a homogeneous self-similar measure on the complex plane increases the correlation dimension (see \ref{dimensiones} for definitions) by a quantitative amount (Theorem \ref{thm:flattening}) and that the Frostman exponent of Bernoulli convolutions tend to 2 as the contraction ratio tends to 1 (Theorem \ref{thm:infty dimension bernoulli}).

The behaviour of the Fourier transform of self-similar measures were also studied in higher dimensions. Let $\mu^{p}_{\lambda,w}$ be the invariant measure associated to the IFS $\{f_j\}_{j=1}^m$, $m\geq 2$ and $f_j = \lambda O x + w_j$, where $\lambda \in (0,1)$, the matrix $O: \mathbb{R}^d \to \mathbb{R}^d$ is an orthogonal matrix and $w_j \in \mathbb{R}^d$ are ``digit'' vectors. Recently, Somolomyak \cite{Sol21} proved that for almost all $\lambda \in (0,1)$, the measure $\mu^{p}_{\lambda,w}$ has power decay at infinity. Indeed the result is more general and is stated for homogeneous  self-affine measures, which are the measures corresponding to an IFS of the form $f_j = Ax + w_j$, where the matrix $A: \mathbb{R}^d\to \mathbb{R}^d$ is a diagonal matrix with the entries $a_{11},\dots, a_{dd}$, with $|a_{ii}|<1$ for all $i=1,\ldots,d$. 
In this work we prove that the Fourier transform of $\mu^p_{\lambda, w}$ has fast power decay outside of a very sparse set of frequencies in the particular case that $O$ is an orthogonal matrix diagonalizable over $\mathbb{R}^d$ (see Proposition \ref{prop:decay_higherdimensions}). 

\section{Fourier decay outside of a sparse set of frequencies}\label{Decay}

Given $p=(p_1,\ldots, p_m)$ a probability vector, i.e a vector in $\mathbb{R}^m$ with $p_i>0$ for all $i=1,\ldots,m$ and $p_1+\ldots +p_m=1$, $w= (w_1,\ldots,w_m)$ a vector sequence in $\mathbb{C}^m$, and $\lambda \in \mathbb{D}$, let $\mu^p_{\lambda,w}$ be the self-similar measure corresponding to the IFS $
\{f_i\}_{i=1}^m$, with probability vector $p$ and where $ f_i= \lambda z+w_i$.
That is, the only Borel probability measure satisfying the relation
\begin{equation*}
\mu^p_{\lambda,w}= \sum_{i=1}^m p_i f_i \mu^p_{\lambda,w}.
\end{equation*}
Here and below, if $\mu$ is a measure on $X$ and $g: X \to Y$ is a map, then $g\mu(A)=\mu(g^{-1}(A))$ is the push-forward measure.\\
On the other hand, $\mu_{\lambda,w}^p$  also can be define as the distribution of the random sum
\begin{equation*}
\sum_{n=1}^{+\infty}  \lambda^n X_n,
\end{equation*}
where $X_n$ are random variables i.i.d with $P(X_n=w_i)= p_i.$

By the definition of $\mu_{\lambda,w}^p$ as a self-similar measure, we can express its Fourier transform as follows
\begin{align*}
\widehat{\mu}_{\lambda,w}^p(\xi)&=\int_{\C} e^{2\pi i \re(z\cdot\overline{\xi})} d\mu_{\lambda,w}^p(z)\\
&= \prod_{n=0}^{\infty} \sum_{j=1}^m p_j \text{exp}(2\pi i \re(\lambda^n w_j\overline{\xi}))\\
&= \prod_{n=0}^{\infty} \Phi(\lambda^n\overline{\xi}),
\end{align*}
where $\Phi(u):= \displaystyle\sum_{j=1}^m p_j \text{exp}(2\pi i \re(w_j u))$.

\begin{lemma}\label{lem:Phi}
The following holds for all $z\in\C$ and $c\in(0,1):$ If $d(\re(z),\Z)>c/2,$ then $|\Phi(z)|<1-\eta(c,p),$ where
\[
\eta(c,p)= p_1+p_2-\sqrt{p_1^2+2p_1p_2\cos(\pi c)+p_2^2}.
\]
\end{lemma}
\begin{proof}
\begin{align*}
|\Phi(z)|&\le |p_1+p_2 e^{2\pi i\re(z)}|+ (1-p_1-p_2)\\
&= |p_1+p_2\cos(2\pi \re(z)) + p_2 i \sin(2\pi \re(z))|+ (1-p_1-p_2)\\
& =\sqrt{p_1^2+p_2^2+2p_1p_2\cos(2\pi \re(z))} + (1-p_1-p_2).
\end{align*}
Now, using that $d(\re(z),\Z)>c/2,$ we obtain that $\cos(2\pi \re(z))<\cos(\pi c)$ and then $|\Phi(z)|\le 1-\eta(c,p)$.
\end{proof}

\begin{proposition} \label{prop:EK0}
Given $\lambda\in\mathbb{D}\setminus\mathbb{R}$ and a probability vector $p=(p_1,\ldots,p_m)$ there is a constant $C=C_{\lambda}>0$ such that the following holds: for each $\varepsilon>0$ small enough (depending continuously on 
$\lambda$) the following holds for all $T$ large enough:  the set of frequencies
$|\xi|\le T$ such that $|\widehat{\mu}_{\lambda,w}^p(\xi)| \ge T^{-\varepsilon}$ can be covered by $C_{\lambda}T^\delta$ squares of side-length 1, where  
\begin{align}
\delta &= \frac{\log\left(\left\lceil \frac{1}{2} \left(1 + \frac{3}{|\lambda|^2}\right)\right\rceil\right)\tilde{\varepsilon}+ h(\tilde{\varepsilon})}{\log(1/|\lambda|)}, \label{eq:def-delta}\\
\tilde{\e} &=\frac{\log(|\lambda|)}{\log\left(1-\eta \left( \frac{|\lambda|^2}{|\lambda|^2+3},p\right)\right)}\varepsilon , \notag
\end{align}
and $h(\tilde{\varepsilon})=-\tilde{\varepsilon}\log(\tilde{\varepsilon})-(1-\tilde{\varepsilon})\log(1-\tilde{\varepsilon})$ is the entropy function.
\end{proposition}
 \begin{proof}
 Choosing $N \in \mathbb{N}$ such that $|\lambda|^{-(N-1)} \leq T \leq|\lambda|^{-N}$ we may assume that $T=|\lambda|^{-N}$. We can write $\xi= \overline{t\lambda^{-N}}$, where $t\in \mathbb{C}$ and $|t|< 1$.
 We have that 
 \begin{align*}
 |\widehat{\mu}_{\lambda,w}^p(\xi)| &\leq \prod_{j=1}^{\infty} |\Phi(\lambda^j\overline{\xi})| \\ 
&= \prod_{j=1}^{\infty} |\Phi(\lambda^j t \lambda^{-N})| \\
&\leq \prod_{j=1}^{N} |\Phi(\lambda^{j-N} t)|\\
&=\prod_{j=0}^{N-1} |\Phi(\lambda^{-j} t)|.
 \end{align*}
 
We denote the distance of $x \in \mathbb{R}$ to the closest integer by $\norm{x}$. Given $\eps>0$, we let $\widetilde{\eps}$ be as in the statement.
  Let
\[
S(N, \tilde{\varepsilon}):=\left\{t\in\C, |t|<1 \colon \|\re(\lambda^{-j} t)\|< \rho \mbox{ for at least } (1-\tilde{\e})N \mbox{ integers } j\in[N]\right\},
\]
where we denote $[N]=\{0,1,\ldots,N-1\}$, and $\rho=\rho(\lambda) =\frac{|\lambda|^2}{2(|\lambda|^2 + 3)}.$

We observe that if $t\notin S(N, \tilde{\varepsilon})$ then, by Lemma \ref{lem:Phi},
\begin{equation}\label{decay}
|\widehat{\mu}_{\lambda,w}^p(\xi)|\le (1-\eta(2|\rho|,p))^{\tilde{\varepsilon} N}=  |\lambda|^{N\e} < T^{-\varepsilon},
\end{equation}
using the definition of $\tilde{\varepsilon}$. We deduce that
\begin{equation} \label{eq:decay-on-SN}
\left\{t\in\C, |t|<1 \colon |\widehat{\mu}_{\lambda,a}^p(\xi)|\ge T^{-\varepsilon}\right\}\subseteq S(N, \tilde{\varepsilon}).
\end{equation}
Hence, in order to prove that $\{\xi\in\C, |\xi|\le T\colon |\widehat{\mu}_{\lambda,a}^p(\xi)|\ge T^{-\varepsilon}\}$ can be covered by a small number of squares of side-length
1, we will estimate the amount and size of squares needed to cover $S(N, \tilde{\varepsilon})$.

For each $t\in\C, |t|<1,$ we define integers $r_j(t)$ and $\varepsilon_j(t)\in [-1/2, 1/2)$ such that
\begin{equation}\label{eq:erdos-kahane}
\re(\lambda^{-j}t)= r_j(t)+\varepsilon_j(t).
\end{equation}
Then $t\in S(N, \tilde{\varepsilon})$ precisely when $|\varepsilon_j(t)|<\rho$ at least $(1-\tilde{\e})N$ times among indices $j\in[N]$.
We will simply write $r_j$ and $\varepsilon_j$ when no confusion arises.

Let $N_1=\lceil (1-\tilde{\e}) N\rceil$. For each $t\in S(N, \tilde{\varepsilon})$, there is a subset $I\subset [N]$ with at least $N_1$ elements such that $|\varepsilon_j|<\rho$ for all $j\in I$. We will estimate the size of 
$S(N,\tilde{\varepsilon})$ by considering each index set $I$ separately, and for this we define
\[
S(I,\tilde{\varepsilon}) = \{t\in\C, |t|<1\colon \|\re(\lambda^{-j}t)\|< \rho \text{ for all } j\in I\}.
\]
We have that $\re(t)=r_0+\varepsilon_0$, so for $|t|<1$ there are at most $3$ choices for $r_0$ and at most $4^3$ choices for $r_1$.

Next, we denote $\lambda= a+ ib$, $a,b \in \mathbb{R}$, $b\neq 0$. Given  $j \in [N]$, we denote \begin{equation}\label{eq:lambda^{-j}t}
 \lambda^{-j}t = c_j + id_j,
\end{equation}
 where $c_j$ and $d_j$ depend on $j$ and $t$. 

By \eqref{eq:erdos-kahane} we have
\begin{align}\label{eq:r_j+1}
\re(\lambda^{-(j+1)}t)&= r_{j+1} + \varepsilon_{j+1}
\end{align}
and, on the other hand 
\begin{align}\label{eq:c_j+1}
c_{j+1} &:= \re(\lambda^{-(j+1)}t)\nonumber \\
&= \re(\lambda^{-1})\re(\lambda^{-j}t)- \im(\lambda^{-1})\im(\lambda^{-j}t) \nonumber \\
&= \frac{a}{|\lambda|^2}(r_j+\varepsilon_j) + \frac{b}{|\lambda|^2}d_j,
\end{align}
where in the last equality we use that $\lambda^{-1}= \dfrac{a-bi}{|\lambda|^2}$.
Using again \eqref{eq:erdos-kahane} and a simple calculation, we obtain
\begin{align}\label{eq:r_j-1}
 r_{j-1}+ \varepsilon_{j-1}= a c_j-b d_j,
\end{align}
and therefore 
\begin{align}\label{eq:d_j}
d_j &= \frac{ac_j-r_{j-1}-\varepsilon_{j-1}}{b}\nonumber\\ &=\frac{1}{b} \left( a(r_j+\varepsilon_j) -r_{j-1} -\varepsilon_{j-1}\right),
\end{align}
where in the last equality we use that $c_j = r_j + \varepsilon_j$.

Now, combining \eqref{eq:r_j+1}, \eqref{eq:c_j+1}, \eqref{eq:d_j} 
\begin{align*}
 \varepsilon_{j+1}= c_{j+1} - r_{j+1}&= \frac{a}{|\lambda|^2}(r_j+\varepsilon_j) + \frac{b}{|\lambda|^2}\left(\frac{1}{b} \left( a(r_j+\varepsilon_j) -r_{j-1} -\varepsilon_{j-1}\right)\right)-r_{j+1}\\&= \frac{2a}{|\lambda|^2}(r_j+\varepsilon_j)-\frac{r_{j-1}+\varepsilon_{j-1}}{|\lambda|^2} - r_{j+1},
\end{align*}
or equivalently,
\[
\frac{2a r_j - r_{j-1}}{|\lambda|^2} - r_{j+1} = \varepsilon_{j+1} + \frac{\varepsilon_{j-1}}{|\lambda|^2} - \frac{2a\varepsilon_j}{|\lambda|^2}. 
\]
Taking absolute value in the previous equality and using that $|\varepsilon_i| \leq 1/2$ for all $i$ and $|a|\leq 1$, we get
\begin{align*}
|r_{j+1} - \frac{2a r_j - r_{j-1}}{|\lambda|^2}| &\leq \frac{1}{2}\left(1 + \frac{1}{|\lambda|^2} + \frac{2|a|}{|\lambda|^2}\right)
\\ &\leq \frac{1}{2} \left( 1 + \frac{3}{|\lambda|^2}\right).
\end{align*}
Therefore, given $r_{j-1},r_j$, we at most can have $\left\lceil \frac{1}{2} \left(1 + \frac{3}{|\lambda|^2}\right)\right\rceil$ choices of $r_{j+1}$.

If $j-1, j, j+1 \in I$, then $|\varepsilon_{j-1}|, |\varepsilon_j|, |\varepsilon_{j+1}| < \rho$ so that $|\varepsilon_{j+1}| + \frac{|\varepsilon_{j-1}|}{|\lambda|^2} + \frac{2|a||\varepsilon_j|}{|\lambda|^2} < 1/2$
and at most one value of $r_{j+1}$ is possible.
Note that
\[
|\{ j\in [N]:j-1,j,j+1\in I\}| \ge N-3|N\setminus N_1|-2 \ge (1-3\widetilde{\varepsilon})N-2.
\]
Thus, the total number of sequences $r_1,\ldots,r_N$ is at most
\[
M_N:= 3\cdot 4^3 \left(\left\lceil \frac{1}{2} \left(1 + \frac{3}{|\lambda|^2}\right)\right\rceil \right)^{3\widetilde{\varepsilon}N+2}.
\]
Invoking \eqref{eq:erdos-kahane} and \eqref{eq:lambda^{-j}t} we have
 \begin{equation}\label{eq:c_N}
 c_N= \re(\lambda^{-N}t) = r_N + \varepsilon_N
 \end{equation}
 with $|c_N - r_N|\leq 1/2$.
 
On the other hand,
 \begin{equation}\label{eq:r_N-1}
 r_{N-1} + \varepsilon_{N-1}=\re(\lambda^{-(N-1)}t)\nonumber = \re(\lambda \lambda^{-Nt}) =  ac_N-bd_N
 \end{equation} 
 and then, using \eqref{eq:r_N-1}, \eqref{eq:c_N} and that $|\varepsilon_j|\leq 1/2$ for all $j$, we obtain 
\begin{align}\label{eq:d_N}
\displaystyle|d_N-\frac{a r_N}{b}- \frac{r_{N-1}}{b}| &\leq \frac{|a|+1}{2|b|}.
\end{align}

\medskip

From \eqref{eq:c_N} and \eqref{eq:d_N} we conclude that, for each pair $(r_{N-1},r_N)$, the complex number $\lambda^{-N}t = c_N + id_N$ belongs to a rectangle of dimensions $\frac{|a|+1}{2|b|}\times 1$. Therefore $t$ is contained in a rectangle of dimensions  $\frac{|a|+1}{2|b|}|\lambda|^N \times |\lambda|^N$, and we obtain that $S(I, \widetilde{\varepsilon})$ can be covered by $M_N$ rectangles of the mentioned size. 

By Stirling's formula, we can estimate $N \choose N_1$ and we see that the number of index sets $I$ is at most $e^{h(\tilde{\varepsilon})N}$ for large enough $N$.
Therefore, $S(N,\tilde{\varepsilon})$ can be covered by 
\[
M_N e^{h(\tilde{\varepsilon})N} 
\]
rectangles of dimensions $\frac{|a|+1}{2|b|}|\lambda|^N \times |\lambda|^N$.
Finally, rescaling, we have that $\{\xi\in\C, |\xi|\le T\colon |\widehat{\nu}_{\lambda,a}^p(\xi)|\ge T^{-\varepsilon}\}$ can be covered by the above number of rectangles of dimensions $\frac{|a|+1}{2|b|}\times 1$. Even more, since each rectangle can be covered by a finite number of squares of side-length 1, we conclude the proof.
\end{proof}

\begin{remark}
Clearly, if $\im(\lambda)=0$, the above proof fails, but the exclusion of $\lambda \in \mathbb{R}$ may appear artificial.
When $m=2$, the case $\lambda \in \mathbb{R}$ reduces back to the family of real Bernoulli convolutions, treated in \cite{MS18}. However, it is interesting to consider the case $m\geq 3$ and the numbers $w_i$ are not collinear.
More precisely, we are going to consider  $\mu_{\lambda,w}^p$ as the self-similar measure corresponding to the IFS $\{f_i=\lambda z+w_i\}_{i=1}^m$ with $\lambda \in (0,1)$, weights $(p_1,\ldots, p_m)$ and the $w_i$  not all collinear. In this case, without lost of generality we can assume that $w_1=0, w_2=1$ and $w_3=i$. 
\end{remark}

\begin{lemma}\label{lem:bound_phi_notcollinear}
The following holds for all $z \in\C$, and $c\in(0,1)$: If $\vertiii{z}>\frac{c}{2}$ then $|\Phi(z)| <1-\eta(c,p)$, where $\eta(c,p)$ is a positive constant and
$\vertiii{\cdot}$ denotes the distance of a complex number to  the lattice $\Z^2$.
\end{lemma}
\begin{proof}
\begin{align*}
| \Phi(z)| &= | p_1 + p_2 \exp(2\pi i \re(z))+ p_3 \exp(2\pi i \re(iz)) + (1-p_1 -p_2-p_3)| \\ &\leq 1- c_1 \max(\|\re(z)\|, \|\re(iz)\|)^2 = 1- c_1 \max(\|\re(z)\|, \|-\im(z)\|)^2 \\& = 1- c_1 \max(\|\re(z)\|, \|\im(z)\|)^2,
\end{align*}
for some constant $c_1>0$ depending on $p$. Here $\|x\|$ denotes the distance from
$x$ to the nearest integer. 
Since
\[
\vertiii{z}^2 = \|\re(z)\|^2 + \|\im(z)\|^2 \leq 2\max(\|\re(z)\|, \|\im(z)\|)^2,
\]
we obtain
\[
|\Phi(z)| \leq 1- \frac{c_1}{2}\vertiii{z}^2 \leq 1-\eta(c,p).
\]
\end{proof}

\begin{proposition} \label{prop:EK0_real}
Given $\lambda\in (0,1)$ and a probability vector $p=(p_1,\ldots,p_m)$, $m\geq 3$, there is a constant $C=C_{\lambda}>0$ such that the following holds: for each $\e>0$ small enough (depending continuously on $\lambda$) the following holds for all $T$ large enough:  the set of frequencies
$|\xi|\le T$ such that $|\widehat{\mu}_{\lambda,w}^p(\xi)| \ge T^{-\e}$ can be covered by $C_{\lambda} T^\delta$ squares of side-length $1$, where
\begin{align*}
\delta &=\frac{4\widetilde{\e} \log(\lceil 2+1/\lambda\rceil) + h(\widetilde{\e})}{\log(1/\lambda)} \\
\widetilde{\e} &= \frac{\varepsilon \log(\lambda)}{\log(1-\eta(\frac{\lambda}{2(\lambda+1)},p))} . 
\end{align*}
\end{proposition}
\begin{proof}
The argument is similar to the argument used in the real case, so we omit the proof.
\end{proof}

\section{A Kaufman's type theorem in two dimensions}
As an application of the previous result we obtain a version of Kaufman's theorem in the complex plane. 

\begin{theorem}\label{Kaufmann on C}
Let $\mu$ be an homogeneous self-similar measure on $\mathbb{C}$ which is not a single atom  and let $F:\mathbb{C} \to \mathbb{C}$ be an analytic function with  $F''$  non zero in a neighborhood of $K:=\supp \mu$. Then there exist $\sigma=\sigma(\mu)>0$ and $C= C(F,\mu)>0$ such that 
\[
|\widehat{F\mu}(\xi)| \leq C |\xi|^{-\sigma}.
\]
\end{theorem}

For the proof we need the following well known result (see for example \cite{FengLau09}).
\begin{proposition}\label{prop:exponent_s}
Let $\mu$ be a self-similar measure on $\mathbb{C}$ which is not a single atom. Then there exist positive constants $C$ and $s$, depending on $\mu$, such that $\mu(B(x,r))\leq Cr^s$, for all $x,r>0.$
\end{proposition}

\begin{proof}[Proof of Theorem \ref{Kaufmann on C}]
Fix  $\xi$ such that $|\xi|\gg1 $ and choose $N\in \mathbb{N}$ such that $1< |\lambda|^N |\xi|^{2/3} \leq |\lambda|^{-1}$. In this case, we also have $\lambda^N \approx |\xi|^{-2/3}$ and $|\xi| |\lambda|^{2N} \approx |\xi|^{-1/3}$.

 Let us decompose $\mu$ in the following way
\[
\mu = \mu_N \ast \nu_N
\]
where $\mu_N = \ast_{n=1}^N \left( \sum_{j=1}^m p_j \delta_{\lambda^n w_j}\right)$ and $\nu_N$ is a rotated scale down copy of $\mu$ by a factor $\lambda^N$.

For the next calculation, we will write $e(z)=e^{-2\pi i \re{z}}$ for simplicity
\begin{align*}
    \widehat{F\mu}(\xi) &= \int_{\mathbb{C}} e^{2\pi i \re(F(v)\bar{\xi})} \, d\mu(v) 
    \\&=  \int_{\mathbb{C}} \int_{\mathbb{C}} e( F(z+w)\bar{\xi}) \, d\mu_N(z) d\nu_N(w) 
    \\&= \int_{\mathbb{C}} \int_{\mathbb{C}} e((F(z)+ F'(z)w + O(|w|^2))\bar{\xi}) \, d\mu_N(z) \,d\nu_N(w) 
    \\&= \int_{\mathbb{C}} \int_{\mathbb{C}} e(F(z)\overline{\xi} + F'(z)w \overline{\xi})e(O(|\xi||w|^2)) \, d\mu_N(z) d\nu_N(w) 
    \\&=\int_{\mathbb{C}} \int_{\mathbb{C}} e(F(z)\bar{\xi} + F'(z)w\bar{\xi}) (1+O(|\xi|| w|^2)) \, d\mu_N(z) \, d\nu_N(w) \\&= \int_{\mathbb{C}} e( F(z)\bar{\xi}) \left( \int_{\mathbb{C}} e(F'(z)w\bar{\xi}) \, d\nu_N(w)\right) \, d\mu_N(z) 
      \\&+ \int_{\mathbb{C}}\int_{\mathbb{C}}  e(F'(z)w\bar{\xi}+ F(z)\bar{\xi}) O(|\xi| |w|^2)\, d\mu_N(z)\, d\nu_N(w) 
    \\&= \int_{\mathbb{C}} e(F(z)\bar{\xi}) \left( \int_{\mathbb{C}} e(F'(z)w\bar{\xi}) \, d\nu_N(w)\right) \, d\mu_N(z) + O(|\xi| |\lambda|^{2N}),
\end{align*}
where in the third equality we replace $F$ by its linear approximation of $F$ ( where, as usual, $O(X)$ denotes a quantity bounded by $CX$ in modulus) and in the  the fifth equality we use that $|e(\delta)-1|=O(\delta)$.

Then, by the assumptions made at the beginning of the proof on $\xi$ and $N$
\begin{align*}
|\widehat{F\mu} (\xi)| &\leq \left| \int_{\mathbb{C}} e(F(z) \xi)\left( \int_{\mathbb{C}} e(F'(z)w \xi) \, d\nu_N(w)\right) \, d\mu_N(z) \right|  + O(|\xi|^{-1/3}) \\&\leq \int_{\mathbb{C}} |\widehat{\nu}_N(F'(z) \xi)| \, d\mu_N(z) + O(|\xi|^{-1/3}) \\ &=\int_{\mathbb{C}} |\widehat{\mu}(\lambda^N F'(z) \xi)| \, d\mu_N(z) +  O(|\xi|^{-1/3}). 
\end{align*}

Consider $T= M |\lambda|^N|\xi|$, where  $M:=\displaystyle\sup_{z \in \, \supp(\mu)}|F'(z)|$ and fix $\e>0$ to be determined later.

Then, by Proposition \ref{prop:EK0}, there is $C=C(\lambda)>0$ such that the set of frequencies $|\xi|\leq T$ for which $|\widehat{\mu}(\xi)| \ge T^{-\e}$ can be covered by $CT^\delta$ squares with sidelength $1$.  Let $I_1,\dots, I_{CT^{\delta}}$ be these squares. Observe that if $\xi\notin\cup_{j=1}^{CT^{\delta}}I_j$ then $|\widehat{\mu}(\xi)|\le T^{-\e}$.

Consider the set 
\[
\Gamma:=\left\{z \in \supp \mu  \colon \lambda^N F'(z) \xi \in\bigcup_{j=1}^{CT^{\delta}}I_j\right\}.
\]
Then
\[
\int_{\mathbb{C}} |\widehat{\mu}(\lambda^N F'(z) \xi)| \, d\mu_N(z) = \int_{\Gamma} + \int_{\Gamma^{\text{c}}} \leq \mu_N(\Gamma) + T^{-\e} \leq \mu_N(\Gamma) + O(|\xi|^{-1/3}).
\]
In order to conclude the proof, we need to prove that $\mu_N(\Gamma) \leq |\xi|^{-\beta}$, for some $\beta>0$.
 
 First, observe that  $\Gamma$ can be rewrite as
  \[\left\{z \in \supp \mu \colon F'(z) \in \bigcup_{j=1}^{CT^{\delta}}J_j\right\}\]
  where $J_j$ are squares of side-length  $|\lambda|^{-N}|\xi|^{-1} \approx |\lambda|^{N/2}$, or which is the same, $\Gamma = \bigcup_{j=1}^{CT^\delta} J'_j$, where $J'_j = (F')^{-1}(J_j) \cap \supp \mu$.
Using that $F''$  non zero in a neighborhood of $K:=\supp \mu$, we have that  
\[
|z_1-z_2| \leq L\,|F'(z_1) - F'(z_2)|  
\]
for all $z_1,z_2 \in J'_j$ and  $L$ is a positive constant depending on $F$.
Then, for each $j$,
\begin{equation}\label{eq: diametroJ'}
\text{diam} \, J'_j   \leq L\, \text{diam}\, J_j \lesssim L|\lambda|^{N/2}.
\end{equation}

On the other hand, since $\nu_N$ is a rotated and scaled down copy of $\mu$ by a factor $\lambda^N$, if the support of $\mu$ is contained in a ball $B_C(0)$, for some $C= C(\lambda, w_1,\ldots, w_m)$, the support of $\nu_N$ is contained in a ball $B_{|\lambda|^NC}(0)$. Then, since $\mu=\mu_N \ast \nu_N$, one have that for any ball $B$
\begin{equation}\label{prop-N}
\mu_N(B) \le \mu(B + B_{\lambda^N C}(0)).
\end{equation}

 Then, invoking \eqref{eq: diametroJ'}, \eqref{eq: diametroJ'} and \eqref{prop-N},  for each $j$ there exists a ball $B_j$ with the same diameter of $J'_j$ such that $J'_j \subseteq B_j$ and $\mu_N(B_j)\leq C|\lambda|^{\frac{N s}{2}}$. 
 
 Therefore
\[
\mu_N(\Gamma) \leq C T^{\delta}|\lambda|^{\frac{N s}{2}} \leq C|\xi|^{\delta/3}|\xi|^{-s/3} \leq C |\xi|^{\frac{(\delta -s)}{3}}.
\]
Choosing $\varepsilon$ small enough such that $\delta < s(\mu)$ we obtain that 
\[
 |\widehat{F\mu}(\xi)| \leq C |\xi|^{\frac{(\delta -s)}{3}} + C|\xi|^{-\varepsilon/3} + C |\xi|^{-1/3} \leq C|\xi|^{-\min\{\frac{(s-\delta)}{3}, \frac{\varepsilon}{3}\}}.
\]
\end{proof}
 
\begin{remark} The above theorem allow us to have uniform explicit power decay for the Fourier transforms of $F\mu_{\lambda}$, the push-forward measure of complex Bernoulli convolutions, even if the measure $\mu_\lambda$ does not have decay at all.  An example are the complex Bernoulli convolutions $\mu_{\lambda}$ with $\lambda = 1/\theta$ and $\theta$ is a complex Pisot number such that $1<|\theta|<\sqrt{2}$. It is known that $|\widehat{\mu_{\lambda}}(\xi)| \nrightarrow 0$ when $|\xi| \to \infty$, see \cite{SX03}. 
\end{remark}

\section{Applications}

\subsection{Improving the $L^2$ dimension under convolution}\label{dimensiones}

We begin by recalling the definition of $L^q$ dimensions. Let $q\in (1,+\infty)$, and set $s_n(\mu,q)= \sum_{Q\in\mathcal{D}_n} \mu(Q)^q$, with $\{\mathcal{D}_n\}_n$ the partition of $\R^d$ into dyadic intervals of length $2^{-n}$. Define
\[
\dim_q(\mu):=\liminf_{n\to+\infty} \frac{\log(s_n(\mu,q))}{(q-1)\log(2^{-n})}.
\]
The $L^2$ dimension of a measure is also known as \emph{correlation dimension}. Note that the Frostman exponent $\dim_\infty$ can also be defined as
\[
\dim_\infty(\mu):=\liminf_{n\to +\infty}\frac{\log \max\{\mu(Q):Q\in\mathcal{D}_n\}}{\log(2^{-n})}.
\]
It is well known that the function $q\mapsto \dim_q(\mu)$ is continuous and non-increasing on $(1,+\infty]$ and that $\dim_q(\mu) \le \hdim(\mu)$ for any $q\in (1,+\infty]$, where $\hdim$ is the lower Hausdorff dimension of a measure, defined as
\[
\hdim(\mu):=\inf\{ \hdim(A):\mu(A)>0\}.
\]
We refer the reader to \cite{FLR02} for the proofs of these facts and further background on dimensions of measures.


The proofs of the results in this section are similar to the ones given in \cite{MS18} for the real case, but we will include them here for the sake of completeness.
\begin{theorem} \label{thm:flattening}
Let $\mu$ be an homogeneous self-similar measure on the complex plane. Given any $\kappa>0$, there is $\sigma=\sigma(\lambda,p,\kappa)>0$ such that the following holds: let $\nu$ be any Borel probability measure with $\dim_2(\nu) \le 2-\kappa$. Then
\[
\dim_2(\mu*\nu) > \dim_2(\nu)+\sigma.
\]
More precisely, one can take $\sigma=2\e$, where $\e=\e(\lambda,p,\kappa)$ is such that the value of $\delta=\delta(\e,\lambda,p)$ given in Proposition \ref{prop:EK0} satisfies
\begin{equation} \label{eq:def-epsilon}
\kappa - 2\e = \delta.
\end{equation}
\end{theorem}
\begin{proof}
First, note that  it is possible to choose $\e<1/2$ such that $\kappa -2\e = \delta$, using continuity arguments.

For any Borel probability measure $\eta$ on $\mathbb{C}$ we have that $\dim_2 \eta = 2-\alpha (\eta)$, where
\[
\alpha(\eta) =\limsup_{T \to \infty}\frac{\log \int_{|\xi|<T} |\widehat{\eta}(\xi)|^2 \, d\xi}{\log T}.
\]
See \cite[Lemma 2.5]{FNW02} for a proof of this result in the real line, but the same argument can be extended to higher dimensions.
Then it is enough to prove that $\alpha(\nu)\geq \kappa$ implies $\alpha (\mu * \nu)< \alpha(\nu) - \sigma$.

Denote $\kappa_0 = \alpha(\nu)\ge \kappa$. For any $\e_0>0$ and taking $T=2^N$ for some $N\in \mathbb{N}$, by definition of $\alpha$,
\[
\int_{ |\xi| \le 2^N}|\widehat{\nu}(\xi)|^2 d\xi \le O_{\e_0}(1) 2^{N(\kappa_0+\e_0)}.
\]
Split the frequencies into two groups
\begin{align*}
E_N &= \{ \xi:  |\xi| \le 2^N, |\widehat{\mu}(\xi)| \le 2^{-\e N}\},\\
F_N &= \{ \xi:  |\xi| \le 2^N, |\widehat{\mu}(\xi)| > 2^{-\e N}\}.
\end{align*}
Then, applying Proposition \ref{prop:EK0}, we have that $F_N$ can be covered by $C_{\lambda} 2^{\delta N}$ squares of side-length $1$ and, in consequence, it has Lebesgue measure bounded by $C_{\lambda} 2^{\delta N}$.

Using all this, we have
\begin{align*}
\int_{|\xi| \le 2^N} |\widehat{\mu*\nu}(\xi)|^2 \,d\xi &= \int_{E_N\cup F_N} |\widehat{\mu}(\xi)|^2|\widehat{\nu}(\xi)|^2\,d\xi \\
&\le \int_{E_N} 2^{-2\e N} |\widehat{\nu}(\xi)|^2 \,d\xi+ \int_{F_N} 1 \,d\xi\\
&\le O_{\e_0}(1) 2^{-2\e N} 2^{(\kappa_0+\e_0)N} + C_\lambda 2^{\delta N}\\
&\le O_{\e_0,\lambda}(1) 2^{(\kappa_0-2\e+\e_0)N},
\end{align*}
using that $\kappa_0\ge \kappa$ and the definition of $\e$ in the last line. Since this holds for all $\e_0>0$, it follows from the definition of $\alpha$ that
\[
\alpha(\mu*\nu) \le \kappa_0 - 2\e,
\]
which gives the claim since $\sigma=2\e$.
\end{proof}

\subsection{Frostman exponent for complex Bernoulli convolutions}

Let $\lambda \in \mathbb{D}$ and $p \in (0,1)$. We denote with $\mu^{p}_\lambda$ the biased complex Bernoulli convolution. This means, $\mu^p_\lambda$ is the self-similar measure associated with the IFS $\{\lambda z -1, \lambda z +1\}$ with probability vector $(p,1-p)$. When $p=1/2$, we just write $\mu_\lambda$ to denote the usual Bernoulli convolution.

\begin{theorem}\label{thm:infty dimension bernoulli} Given $p_0<1/2$ there exists a constant $C=C(p_0)$ such that 
\[
\dim_{\infty}(\mu^p_\lambda) \geq 2 - C(1-|\lambda|)\log(1/(1-|\lambda|))
\]
for all $p_0 \leq p\leq 1-p_0$.
\end{theorem}

\begin{proof}
Fix $\lambda \in \mathbb{D}$, with modulus close to $1$. We define $N= N(\lambda)$ the smallest integer such that $|\lambda|^N <1/\sqrt{2}$.

Then,
\[
\frac{|\lambda|}{\sqrt{2}}\le |\lambda|^N<\frac{1}{\sqrt{2}}.
\]
In particular, assuming $|\lambda|>1/\sqrt{2}$, we see that $|\lambda|^N\in (1/2,1/\sqrt{2})$.

Fix $\kappa\in (0,1)$, and suppose that $\dim_2(\mu_\lambda^p)\le 2-\kappa$.
We have the following decomposition. Let us write $S_{\lambda}(x)=\lambda x$ for the map that scales by $\lambda$, and recall that
\begin{equation} \label{eq:convolution-structure}
\mu_\lambda^p = \mu_{\lambda^N}^p * S_\lambda \mu_{\lambda^N}^p * \cdots * S_{\lambda^{N-1}} \mu_{\lambda^N}^p.
\end{equation}
This is a well-known fact that can be seen from expressing $\mu^p_\lambda$ as an infinite convolution.

Since the associated IFS satisfies the open set condition,
\begin{equation} \label{eq:corr-dim-OSC}
\dim_2(\mu_{\lambda^N}^p) = \frac{\log(p^2+(1-p)^2)}{\log(|\lambda|^N)}.
\end{equation}
In particular, we have that $\dim_2(\mu_{\lambda^N}^p)\ge 0$.
Now, using that $\dim_2(\mu_\lambda^p)\le 2-\kappa$ and \eqref{eq:convolution-structure},
we get $\dim_2(\mu_{\lambda^N}^p)\le 2-\kappa$. By Theorem \ref{thm:flattening}, there is $\sigma=\sigma(\lambda, p, \kappa)>0$ such that
\[
\dim_2(\mu_{\lambda^N}^p *S_\lambda \mu_{\lambda^N}^p)\ge \sigma.
\]
Proceeding inductively according to \eqref{eq:convolution-structure}, after $N-1$ steps we obtain that if $\dim_2(\mu_\lambda^p)\le 2-\kappa$, then
\[
\dim_2(\mu_\lambda^p) \ge (N-1)\sigma.
\]
It follows that if $\kappa$ is such that $\sigma=\sigma(\lambda, p,\kappa)=1/(N-1)$, then
\[
\dim_2(\mu_\lambda^p)\ge 2-\kappa.
\]
Thus, we just need to estimate such $\kappa$.  By \eqref{eq:def-epsilon}, we have that $\kappa = \delta+\sigma$, where $\delta=\delta(\sigma/2)$ is given by \eqref{eq:def-delta}. Note that $\tilde{\e}=C(\lambda^N,p)\sigma$, where $C>0$ depends continuously on $\lambda^N$ and $p$. In what follows, $C_j$ will denote a postive constant depending only of $p_0$.  Since $|\lambda|^N\in (1/2,1/\sqrt{2})$ and $p\in [p_0,1-p_0]$, a calculation using \eqref{eq:def-delta} shows that there is a constant $C_1$ such that $\delta\le C_1 \sigma\log(1/\sigma)$ provided $\sigma$ is small enough (which we may assume). 
We deduce that
\begin{equation} \label{eq:final-lower-bound}
\dim_2(\mu_\lambda^p)\ge 2-\kappa \ge 2-\sigma-C_1\sigma\log(1/\sigma) \ge 2- C_2\sigma\log(1/\sigma)
\end{equation}
if $\sigma$ is small enough. On the other hand, since $|\lambda|^{1/\sigma}=|\lambda|^{N-1}<|\lambda|^{-1}/\sqrt{2}<2/3$ (say), we have that $\sigma \le \log(1/|\lambda|)/\log(3/2)$. Finally, using that $\log(1/|\lambda|) \le 2(1-|\lambda|)$ for $1-|\lambda|$ small, we deduce that
\[
\sigma \le C_4 (1-|\lambda|).
\]
Together with \eqref{eq:final-lower-bound}, this yields
\[
\dim_2(\mu_\lambda^p) \ge  2 - C_5 (1-|\lambda|) \log(1/(1-|\lambda|)).
\]
Since we have the decomposition
\[
 \mu_\lambda^p = \mu_{\lambda^2}^p * S_\lambda\mu_{\lambda^2}^p,
\]
and scalings do not change $L^2$ dimension, using Young's Lemma (see Lemma 5.2 in \cite{MS18}), we can conclude that
\[
\dim_\infty(\mu_\lambda^p) \ge 2 - C_5 (1-|\lambda|^2) \log(1/(1-|\lambda|^2)) \ge 2 - C_6(1-|\lambda|)\log(1/(1-|\lambda|)).
\]
\end{proof}

\begin{corollary} \label{cor:unbiased}
There is an absolute constant $C>0$ such that
\[
\dim_\infty(\mu_\lambda) \ge 2 - C(1-|\lambda|)^2 \log(1/(1-|\lambda|)).
\]
\end{corollary}
\begin{proof}
Again, fix $\lambda \in \mathbb{D}$ with $|\lambda|$ close to $1$ and let $N=N(\lambda)$ the smallest integer such that $|\lambda|^N < 1/\sqrt{2}$ and then, $|\lambda|^N > |\lambda|/\sqrt{2}$.

Since the associated IFS satisfy the open set condition,
\[
\dim_2(\mu_{\lambda^N}) = \frac{\log(1/2)}{\log |\lambda|^N} \ge \frac{\log (1/2)}{\log( |\lambda|/\sqrt{2})} = 1 -  \frac{\log(1/ |\lambda|)}{\log(2/ |\lambda|)}.
\]

Proceeding as in the proof of the Theorem \ref{thm:infty dimension bernoulli}, we obtain that if $\dim_2(\mu_{\lambda})\le 2- \kappa$ then there exists
$\sigma=\sigma(\kappa, \lambda)>0$ such that
\begin{align*}
\dim_2(\mu_{\lambda})&\ge \dim_2(\mu_{ \lambda^N})+ (N-1)\sigma\\
&\ge 1-\frac{\log(1/ |\lambda|)}{\log(2/ |\lambda|)}+ (N-1)\sigma.
\end{align*}
Now if $\kappa$ is such that $\sigma=\frac{\log(1/ |\lambda|)}{\log(2/ |\lambda|)}\frac{1}{N-1},$ then
\[
\dim_2(\mu_{\lambda})\ge 2-\kappa.
\]
Then we want to estimate such $\kappa$. Proceeding as in the proof of the above theorem, we get
\[
\dim_2(\mu_{\lambda})\ge 2-\kappa\ge 2-C_1\sigma\log(1/\sigma).
\]
On the other hand, using that $|\lambda|/\sqrt{2}< |\lambda|^N<1/{\sqrt{2}}$ and that $\log(1/ |\lambda|) \le 2(1- |\lambda|)$ for $1- |\lambda|$ small, we obtain  $\frac{1}{N-1}\le C_2(1- |\lambda|)$ and then
 \[
 \sigma\le C_3(1- |\lambda|)^2.
 \]

Thus
 \[
\dim_2(\mu_{\lambda}) \ge 2 - C_4(1-|\lambda|)^2 \log(1/(1-|\lambda|)).
\]

Using again Young's Lemma as in the proof of Theorem \ref{thm:infty dimension bernoulli} we finish the proof.

\end{proof}

\section{A particular extension to higher dimensions}\label{Rn}

Let $O: \mathbb{R}^d \to \mathbb{R}^d$ an orthogonal matrix transformation diagonalizable over $\mathbb{R}^d$. Given $(p_1,\ldots,p_m)$ a probability vector, $w=\{w_i\}_{i=1}^m$ a sequence of digit vectors in $\mathbb{R}^d$, with $m\geq 3$ and  $\lambda \in (0,1)$, let $\mu^p_{\lambda,w}$ the self-similar measure associated to the IFS $\{f_i\}_{i=1}^m$, where $f_i = \lambda Ox + w_i$.

Since $O \in O(n)$ is diagonalizable over $\mathbb{R}^n$, we can assume without lost of generality that $O$ is the identity matrix. In fact, if we iterate the IFS (replacing $f_i$ by $(f_if_j)_{i,j=1}^m$), the matrix $O$ is replaced by $O^2$ which has all its eigenvalues equal to 1 and then, in some autonormal basis, it can be written as the identity matrix. 
Then, after an affine change of coordinates we can always assume that $w_1=(0,\ldots,0)$ and $w_2=(1,\ldots,1)$.
By definition of the Fourier transform and the condition of being a self-similar measure we can write 
\[
\widehat{\mu}^p_{\lambda,w}(\xi) = \prod_{n=0}^{\infty} \Phi(\lambda ^n\xi), 
\]
where 
\[
\Phi(y)=\sum_{j=1}^m p_j e^{-2\pi i \langle y, a_j \rangle\! \ }, \, \, \text{for all} \, \, y \in \mathbb{R}^d.
\]
 
\begin{lemma}\label{lem:phi_higherdimensions}
The following holds for all $y\in\mathbb{R}^d$, $y=(y_1,\ldots, y_d)$ and $c\in(0,1)$: If $||y_1 +\ldots + y_n||>\frac{c}{2}$ then $\Phi(y) <1-\eta(c,p)$, where $\eta(c,p)$ is a positive constant and  
$||\cdot||$ denotes the distance of a real number to the closest integer.
\end{lemma}

\begin{proposition} \label{prop:decay_higherdimensions}
Given $\lambda\in (0,1)$ and a probability vector $p=(p_1,\ldots,p_m)$, $m\geq 3$, there is a constant $C=C_{\lambda}>0$ such that the following holds: for each $\e>0$ small enough (depending continuously on $\lambda$) the following holds for all $T$ large enough:  the set of frequencies
$\{\|\xi\|_{\infty}\le T \colon |\widehat{\mu}_{\lambda,w}^p(\xi)| \ge T^{-\e}\}$ can be covered by $C_{\lambda} T^{\delta}$ squares of side-length $1$, where 
\begin{align*}
\delta &=\frac{\log(\lceil 1+1/\lambda\rceil)\tilde{\varepsilon}+h(\tilde{\varepsilon})}{\log(1/{\lambda})} \\
\widetilde{\e} &=\frac{\log(\lambda)}{\log(1-\eta(\frac{\lambda}{\lambda+1}, p))} \varepsilon.
\end{align*}
\end{proposition}

The proof of Lemma \ref{lem:phi_higherdimensions} is analogous to the proof of Lemma \ref{lem:Phi} and the proof of Proposition \ref{prop:decay_higherdimensions} is similar to the proof of the real case in \cite{MS18}.

\section{Acknowledgements}
The authors would like to thank Pablo Shmerkin for bring to us this problem and for his extremely valuable and helpful comments and suggestions.


\begin{thebibliography}{10}

\bibitem{DavenportEtAl63}
H.~Davenport, P.~Erd\H{o}s, and W.~J. LeVeque.
\newblock On {W}eyl's criterion for uniform distribution.
\newblock {\em Michigan Math. J.}, 10:311--314, 1963.

\bibitem{FLR02}
Ai-Hua Fan, Ka-Sing Lau, and Hui Rao.
\newblock Relationships between different dimensions of a measure.
\newblock {\em Monatsh. Math.}, 135(3):191--201, 2002.

\bibitem{FengLau09}
De-Jun Feng and Ka-Sing Lau.
\newblock Multifractal formalism for self-similar measures with weak separation
  condition.
\newblock {\em J. Math. Pures Appl. (9)}, 92(4):407--428, 2009.

\bibitem{FNW02}
De-Jun Feng, Nhu~T. Nguyen, and Tonghui Wang.
\newblock Convolutions of equicontractive self-similar measures on the line.
\newblock {\em Illinois J. Math.}, 46(4):1339--1351, 2002.

\bibitem{Kaufman84}
Robert Kaufman.
\newblock On {B}ernoulli convolutions.
\newblock In {\em Conference in modern analysis and probability ({N}ew {H}aven,
  {C}onn., 1982)}, volume~26 of {\em Contemp. Math.}, pages 217--222. Amer.
  Math. Soc., Providence, RI, 1984.

\bibitem{Mitsis02}
Themis Mitsis.
\newblock A Stein-Tomas restriction theorem for general measures.
\newblock {\em  Publ. Math. Debrecen.}, 60 (2002), no. 1-2, 89--99.

\bibitem{Mockenhaupt00}
Gerd Mockenhaupt.
\newblock Salem sets and restriction properties of {F}ourier transforms.
\newblock {\em Geom. Funct. Anal.}, 10(6):1579--1587, 2000.

\bibitem{MS18}
Carolina A. Mosquera and Pablo Shmerkin 
\newblock Self-similar measures: asymptotic bounds for the dimension
and Fourier decay or smooth images.
\newblock {\em Ann. Acad. Sci. Fenn. Math.},  43(2):823--834, 2018.

\bibitem{RS20}
Eino Rossi and Pablo Shmerkin
\newblock On measures that improve $L^q$ dimension under convolution. 
\newblock {\em Rev. Mat. Iberoam.},  36 (7): 2217--2236, 2020.

\bibitem{QR03}
 Martine Queff\'{e}lec and Olivier Ramar\'{e}
 \newblock Analyse de {F}ourier des fractions continues \`a quotients restreints.
\newblock {\em Enseign. Math. (2), 49(3-4):335–356, 2003.} 
 

\bibitem{SS16}
Pablo Shmerkin and Boris Solomyak
\newblock Absolute continuity of complex Bernoulli convolutions.
\newblock {\em Math. Proc. Cambridge Philos. Soc.}, 161(3):435--453, 2016.

\bibitem{Shmerkin16}
Pablo Shmerkin.
\newblock On {F}urstenberg's intersection conjecture, self-similar measures,
  and the ${L}^q$ norms of convolutions.
\newblock Preprint, arXiv:1609.07802, 2016.

\bibitem{Tsujii15}
Masato Tsujii.
\newblock On the {F}ourier transforms of self-similar measures.
\newblock {\em Dyn. Syst.}, 30(4):468--484, 2015.

\bibitem{Sol21}
Boris Solomyak.
\newblock Fourier decay for homegenous self-affine measures.
\newblock Preprint, arxiv 2105.08129, 2021.

\bibitem{SX03}
Boris Solomyak and Hui Xu
\newblock On the ``Mandelbrot set" for a pair of linear maps and complex Bernoulli convolutions.
\newblock {\em Nonlinearity.}, 16(5):1733--1749, 2003.

\end{thebibliography}
\end{document}